\documentclass[12pt,amssymb]{article} 

\usepackage{mathrsfs}
\usepackage{amssymb}
\usepackage[dvips]{graphicx}

\newcommand{\newsection}[1]{
\addtocounter{section}{1} \setcounter{equation}{0}
\setcounter{subsection}{0} \addcontentsline{toc}{section}{\protect
\numberline{\arabic{section}}{{\rm #1}}} \vglue .6cm \pagebreak[3]
\noindent{ \bf  \thesection. #1}\nopagebreak[4]\par\vskip .3cm}
\newcommand{\newsubsection}[1]{
\addtocounter{subsection}{1}\setcounter{subsubsection}{0}
\addcontentsline{toc}{subsection}{\protect
\numberline{\arabic{section}.\arabic{subsection}}{#1}} \vglue .4cm
\pagebreak[3] \noindent{\it \thesubsection.
#1}\nopagebreak[4]\par\vskip .3cm}

\makeatletter
\newcommand{\seclabel}[1]{%
  \@bsphack
  \protected@write\@auxout{}%
     {\string\newlabel{#1}{{\thesection}{\thepage}}}
  \@esphack
  }
\newcommand{\subseclabel}[1]{%
  \@bsphack
  \protected@write\@auxout{}%
     {\string\newlabel{#1}{{\thesubsection}{\thepage}}}
  \@esphack
  }
\newcommand{\tablabel}[1]{%
  \@bsphack
  \protected@write\@auxout{}%
     {\string\newlabel{#1}{{\arabic{tabnum}}{\thepage}}}
  \@esphack
  }
\makeatother

\renewcommand{\theequation}{\thesection.\arabic{equation}}

\newlength{\extraspace}
\newlength{\extraspaces}
\newcounter{dummy}
\newcommand{\bc}{\begin{center}}
\newcommand{\ec}{\end{center}}
\newcommand{\be}{\begin{equation}
\addtolength{\abovedisplayskip}{\extraspaces}
\addtolength{\belowdisplayskip}{\extraspaces}
\addtolength{\abovedisplayshortskip}{\extraspace}
\addtolength{\belowdisplayshortskip}{\extraspace}}
\newcommand{\ee}{\end{equation}}

\newcommand{\ba}{\begin{eqnarray}
\addtolength{\abovedisplayskip}{\extraspaces}
\addtolength{\belowdisplayskip}{\extraspaces}
\addtolength{\abovedisplayshortskip}{\extraspace}
\addtolength{\belowdisplayshortskip}{\extraspace}}
\newcommand{\ea}{\end{eqnarray}}

\newcommand{\ban}{\begin{eqnarray*}
\addtolength{\abovedisplayskip}{\extraspaces}
\addtolength{\belowdisplayskip}{\extraspaces}
\addtolength{\abovedisplayshortskip}{\extraspace}
\addtolength{\belowdisplayshortskip}{\extraspace}}
\newcommand{\ean}{\end{eqnarray*}}
\newcommand{\baa}{
\addtocounter{equation}{1} \setcounter{dummy}{\value{equation}}
\setcounter{equation}{0}
\renewcommand{\theequation}{\thesection.\arabic{dummy}\alph{equation}}
\begin{eqnarray}
\addtolength{\abovedisplayskip}{\extraspaces}
\addtolength{\belowdisplayskip}{\extraspaces}
\addtolength{\abovedisplayshortskip}{\extraspace}
\addtolength{\belowdisplayshortskip}{\extraspace}}
\newcommand{\eaa}{
\end{eqnarray}
\setcounter{equation}{\value{dummy}}
\renewcommand{\theequation}{\thesection.\arabic{equation}}}

\input epsf

\newcounter{fignum}
\newcounter{tabel}

\newcounter{tabnum}
\newcommand{\cO}{{\cal O}}
\newcommand{\Ext}{{\rm Ext}}
\newcommand{\Hom}{{\rm Hom}}

\begin{document}

%
%

\begin{flushright}
October 2015\\
\end{flushright}
\vspace{2cm}

\thispagestyle{empty}

\begin{center}
{\Large\bf ADE Transform
 \\[13mm] }

{\sc Ron Donagi \& Martijn Wijnholt}\\[2.5mm]
{\it Department of Mathematics, University of Pennsylvania \\
Philadelphia, PA 19104-6395, USA}\\[15mm]

Abstract:
\end{center}

There is a beautiful correspondence between configurations of lines on a rational surface
and tautological bundles over that surface. We extend this correspondence to families,
by means of a generalized Fourier-Mukai transform that relates spectral data to
bundles over a rational surface fibration.

\newpage

\renewcommand{\Large}{\normalsize}

\tableofcontents

\newpage

\newsection{General comments}

It has long been appreciated that there are remarkable relations between del Pezzo surfaces
and exceptional Lie groups. One way this manifests itself is through a correspondence
between configurations of lines on a del Pezzo surface, and tautological $\widetilde{E}_n$ bundles over
the surface, where $\widetilde{E}_n = E_n \times_{{\cal Z}} {\bf C}^*$ and ${\cal Z}$ is
the center of $E_n$.

For example, the homology classes of the twenty-seven lines $l_i$ on a cubic surface
are associated to the twenty-seven dimensional
representation of $E_6$. The corresponding bundle
\be
V\ =\ \bigoplus_{i=1}^{27} \cO(l_i)
\ee
is an $\widetilde{E}_6$ vector bundle
of rank twenty-seven on the cubic surface,
where $\widetilde{E}_6 = E_6 \times_{{\bf Z}_3} {\bf C}^*$.

When restricted to a fixed elliptic curve embedded as the anti-canonical divisor
of the del Pezzo, and twisted to get a degree zero bundle,
the structure group of the bundle reduces to $E_n$.
By varying the surface,
this correspondence yields all possible semi-stable
$E_n$ bundles over an elliptic curve
\cite{Looijenga:ell_sing,Looijenga:ell_inv,Friedman:1997ih}.

This construction has played an important role in heterotic/$F$-theory duality
\cite{Morrison:1996na,Friedman:1997ih,Clingher:2003ui}.
Moreover the correspondence is not restricted to exceptional groups,
and can be extended to other Lie groups by considering more general rational surfaces \cite{LZ1,LZ2},
for essentially all one needs in the construction
is a rational surface whose homology lattice contains an ADE sublattice.

Our purpose in this work is to extend this correspondence to families of rational surfaces.
In order to accomplish this, we will introduce an integral transform for
rational surface fibrations. It plays a role somewhat analogous to
that of the Fourier-Mukai transform for bundles over elliptically fibered
manifolds.

For convenience let us briefly recall some aspects of the Fourier-Mukai transform.
Given an elliptic fibration
\be
\pi_Z: Z\ \to\ B,
\ee
the Fourier-Mukai transform \cite{Huyb_FM,BBR_FM}  takes a sheaf
$\cal{F}$ on $Z$, or a
complex of sheaves, to another such complex:
\be\label{FM_transform_def}
{\bf FM}({\cal F}) \ = \  Rp_{1*}(p_2^*{\cal F} \otimes {\cal P}),
\ee
where ${\cal P}$ is the (relative) Poincare sheaf on the fiber product
$Z \times_B Z$, and $p_1$ and $p_2$ are the projections to the first and second factor
of $Z \times_B Z$ respectively. A particularly interesting case is when the sheaf
${\cal F}$ is a
stable vector bundle on $Z$ and the transform ${\bf FM}({\cal F})$
is a line bundle $L$ on
a spectral cover $C \subset Z$. Conversely, stable vector bundles
${\cal F}$ on $Z$
can be constructed as the Fourier-Mukai transform of spectral data
$(C,L)$, consisting of a line bundle $L$ on a spectral cover $C \subset Z$.
This powerful algebraic tool is useful also for understanding many aspects of heterotic string
compactifications \cite{Friedman:1997ih,Donagi:1998vx}
and $F$-theory \cite{DonagiCovers,Donagi:2008ca}.

We extend this correspondence in the following way.
Given an elliptic fibration $\pi_Z: Z \to B$
and a spectral cover $C \subset Z$, we construct a  fibration
\be
\pi_Y: Y_C\ \to\ B
\ee
whose fibers are (essentially) Leung's rational surfaces, with `boundary' $Z
\subset Y_C$. Here $Y_C$ depends on the spectral cover $C$, but not on the
spectral sheaf $L$.
Our transform ${\bf T}$ then takes a line bundle $L$ supported on the spectral cover
$C
\subset Z$ to a sheaf (or sometimes a complex) ${\bf T}(L)$ supported on $Y_C$.
Furthermore, its restriction to $Z \subset Y_C$ is the usual Fourier-Mukai transform
${\bf FM}(L)$ of $L$.

As for all integral transforms, our transform is given by a variant of
formula (\ref{FM_transform_def}). In our description the analogue of the Poincar\'e sheaf
is asymmetric: it lives on the fiber product
\be
\widehat{Y}_C\ =\ Y_C  \times_B C,
\ee
which is a rational surface fibration over $C$.

The points in the fiber of $C\to B$ parametrize lines in the corresponding rational surface.
As a result there is a `universal line', which is a divisor
$D \subset \widehat{Y}_C$. The
Poincar\'e sheaf is given, in (\ref{P_sheaf_def}), as a twisted version of $\cO(D)$. Its
restriction to $Z \times_B C$ is also the restriction of the usual
Poincar\'e sheaf, on $Z \times_B Z$. The identification of ${\bf FM}(L)$ with the
restriction of ${\bf T}(L)$ follows from this. We work mostly with the $A_n$ case.
A few remarks are made about $D_n$ and $E_n$, but detailed analysis of those
cases is left for a future work.

This construction is clearly `right' at points of $B$ over which the cover
$C \to B$ is unramified. It is more interesting to figure out what happens
over the branch locus. In this direction we give two results. In
(\ref{P_loc}) we show that near a generic branch point our ${\bf T}(L)$ is still
locally free, i.e. it is an honest vector bundle.
We further show in (\ref{P_ext}) that at a generic branch point,
${\bf T}(L)$ is in fact the regular representative, as defined in
\cite{FriedmanMorgan}.

The construction has some rather interesting string-theoretic implications,
which we develop in a separate paper \cite{Bundles_Forms}.
The basic point is that the traditional geometric
formulation of $M$-theory and $F$-theory
compactifications is not able to capture certain non-abelian Hodge
theoretic structures that are present in dual formulations. The ADE transform
shows that such structures can in fact be described, provided we extend
the traditional formulation of an $M$-theory or $F$-theory
compactification to one with bundles or sheaves.
Although at first sight this might seem like a radical proposal,
as a result we can resolve a number of previously problematic issues involving
missing branches of the moduli space, Yukawa couplings and stability conditions.

\newpage

\newsection{ADE surfaces and canonical ADE bundles}

\seclabel{Correspondence}

\newsubsection{Lines and bundles}

In this section we give a brief review of configurations of lines on a rational surface,
and bundles over that surface. In particular for $G$ a simple Lie group
of type ADE, we review
the relation between lines on $G$-surface $S_G$ as defined in \cite{LZ1,LZ2},
canonical $G$-bundles on $S_G$, and spectral data on a `boundary' elliptic curve $E\subset S_G$.

The idea is as follows. Suppose we are given a surface $S$
for which the group of line bundles ${\rm Pic}(S)$ is discrete, i.e.
${\rm Pic}(S) \cong H^2(S,{\bf Z})$, like for a rational surface.
Suppose we are also given a lattice $N$, and define
\be T\ =\ N \otimes {\bf C}^* \ee
Then up to isomorphism, $T$-bundles on $S$ are classified by
\be
\Hom(N^\vee,{\rm Pic}(S))
\ee
where
\be
N^\vee \ = \ \Hom(N,{\bf Z}) \ = \ \Hom(T,{\bf C}^*)
\ee
is the dual of $N$.

Now suppose that $H^2(S,{\bf Z})$ has a sublattice which is
isomorphic to a root lattice $\Lambda_{rt}$ for a
simple ADE Lie group $G$.
We will take $N = \Lambda_{wt}$ so that $N^\vee =
\Lambda_{rt}$. Then, we get a canonical element of $\Hom(N^\vee,
{\rm Pic}(S))$, namely the inclusion of $\Lambda_{rt}$
into $H^2(S,{\bf Z})$,
and therefore (up to isomorphism) we get a canonical $T$-bundle on $S$.
More precisely, the map (and hence the bundle)
is canonical up to the action by the Weyl group. But $T$-bundles
related by an action of the Weyl group determine the same
$G$-bundle, where we identify $T$ with a maximal torus of $G$. So
we get a canonical $G$-bundle on $S$.

Although this correspondence is quite general, we are going to consider some
particular rational surfaces which one may informally think of as a compactification
of ADE type ALE spaces by an elliptic curve $E$ at infinity. The elliptic
curve will arise as an anti-canonical divisor, so the pair
$(S,E)$ will be a log Calabi-Yau manifold, a K\"ahler surface whose log canonical bundle
\be
K_{(S,E)} \ \equiv \ K_S + E
\ee
is trivial. A log Calabi-Yau $(S,E)$ admits a holomorphic volume form with a simple pole
along $E$, which is interpreted as a holomorphic version of the boundary of $S$.
Such manifolds give a large and mostly overlooked class
of string compactifications where we get an interesting interplay
between bulk and boundary. The boundary behaves as a kind of holographic
screen in the sense that holomorphic constructions over $S$
can be translated to constructions over $E$.

Thus we
can restrict our canonical $G$-bundle on $S$ to get a $G$-bundle on
$E$. If furthermore $[E]\in \Lambda_{rt}^\perp \subset {\rm Pic}(S)$ in
the sense that $\alpha(E) = 0$ for any $\alpha \in \Lambda_{rt}$,
then we get a flat $G$-bundle on $E$. As discussed in
\cite{LeungADE,LZ1,LZ2}, all flat $G$-bundles on $E$ may be
recovered in this way, and we can pick an essentially
unique rational surface $S_G$ such that the moduli space of the flat
$G$-bundle on $E$ equals the complex structure moduli space of
pairs $(S,E')$ with an isomorphism $E'\to E$. The surface $S_G$ is rational and can be
constructed very explicitly by a sequence of blow-ups. We will refer to $S_G$ as
a rational
surface of type $G$. For $G = E_k$ these are just the del Pezzo
surfaces, with one extra blow-up as defined in \cite{LZ1,LZ2}.
In string theory applications one would typically add a few more blow-ups,
in order to embed the correspondence in heterotic/$F$-theory duality.

The ADE lattice can also be obtained from the homology classes of lines
living on the rational surface. Then
we can try to construct
associated vector bundles ${\cal V}_\rho$ for each representation $\rho$ of
$G$ by considering associated configurations of lines.
By restriction, they would yield associated bundles on $E$. The homology
classes of the lines however generate a lattice that is slightly larger
than the ADE root lattice. As a result, one typically has to deal with vector
bundles that have a slightly larger structure group of the form
$G_{ADE} \times_{\cal Z} {\bf C}^*$, where ${\cal Z}$ is the center of $G_{ADE}$.

Let us discuss the $D_n$ and $A_{n-1}$ cases in more detail. We take a Hirzebruch surface ${\bf F}^1$ and
blow up in $n$ points to get the $D_n$ surface $S_{D_n}$. The homology lattice is spanned by
$b,f$ and the exceptional curves $l_i$, and the canonical divisor is given by
$K = -2b-3f +\sum l_i$. To get the $A_{n-1}$ surface we can contract $b$
to get $S_{A_{n-1}}$ (or equivalently start with ${\bf P}^2$ and blow up
$n$ points on that). However in practice we will keep $b$ around and also denote it as $l_0$.
The root lattice is given by
\be
\Lambda_{rt}\ =\ \{ x \in {\rm Pic}(S)|x\cdot K = x\cdot f = x\cdot b = 0 \}
\ee
Since $x\cdot b = 0$, this descends to $S_{A_{n-1}}$. The lattice $\Lambda_{rt}$ is a root lattice of
$A_{n-1}$ type, and the simple roots
may be taken as
\be
\alpha_1 = l_1 - l_2, \ \ldots\ , \alpha_{n-1} = l_{n-1}-l_n
\ee
The elliptic curve $E$ is obtained by embedding into ${\bf P}^2$ using the linear
system $|3p_0|$, where $p_0$ is the identity on $E$. Then we blow up ${\bf P}^2$
at $p_0$ and generic points $p_i$ on $E$ to get $l_i$. For $SU(n)$ we also want to impose
that $\sum_{i=1}^n p_i = np_0$ in the group law on $E$. Then the tautological $A_{n-1}$ bundle
on $S$ for the adjoint representation is given by
\be
\cO^{\oplus n-1} \oplus \sum_{\alpha \in \Phi} \cO_S(\alpha)
\ee
The adjoint representation is not faithful, and the bundle above may be equally regarded
as an $SU(n)$ bundle as well as a $U(n) = SU(n)\times_{{\bf Z}_n} U(1)$ bundle.
The tautological bundle associated to the fundamental representation on the other hand
is given by
\be
{\cal W} = \sum_{i=1}^n\cO_S(l_i)
\ee
which is merely a $U(n)$ bundle.
It is convenient to twist this bundle by tensoring with $\cO_S(-l_0)$:
\be
\widetilde{\cal W}\ =\ \sum_{i=1}^n\cO_S(l_i-l_0)
\ee
Although still a $U(n)$ bundle on our rational surface $S$, we have $(l_i-l_0)\cdot E = 0$
so each of the factors restricts to a degree zero bundle on $E$. The restricted bundle
\be
{\cal V}\ =\ \cO_E(p_1-p_0) \oplus \ldots \oplus \cO_E(p_n-p_0)
\ee
on $E$ does naturally allow a reduction of the structure group from $U(n)$ to $SU(n)$.

For $D_n$, the natural faithful representation to consider is the $2n$-dimensional vector
representation. The associated vector bundle can be constructed from the lines as follows:
\be
{\cal W}\ =\ \bigoplus_{i=1}^n\cO_S(l_i) \oplus \cO_S(f - l_i)
\ee
Its first Chern class is non-vanishing, so the structure group is actually
a $U(1)$ extension of $SO(2n)$. Considering the twisted version
\be
\widetilde {\cal W}\ =\ \bigoplus_{i=1}^n\, \cO_S(l_i-l_0) \oplus \cO_S(f - l_i-l_0)
\ee
we get the following bundle when restricted to $E$:
\be
{\cal V}\ =\ \bigoplus_{i=1}^n \cO_E(p_i-p_0) \oplus \cO_E(-p_i+p_0)
\ee
The bundle ${\cal V}$ does admit the structure of a flat $SO(2n)$ bundle.

\newsubsection{Singular surfaces and the regular representative}

\subseclabel{Reg_representative}

The correspondence reviewed above works for generic, smooth rational surfaces.
Now we would like to ask what happens to this correspondence when the surface is allowed to
degenerate. This is a very common situation and the singularities often have interesting
physics associated to them.
Remarkably, the construction above extends to singular rational surfaces
with singularities of ADE type. Physically such singularities are usually associated to enhanced
gauge symmetries. Moreover such singularities are generically unavoidable when
we study families of rational surfaces, so understanding what happens will
be an important aspect of our construction.

Let us first consider this from the point of view of the boundary elliptic curve $E$.
For simplicity we consider the rank two vector bundle  $\cO_E(p_1)\oplus \cO_E(p_2)$.
Suppose that we take the limit as $p_2 \to p_1$. By the classification of bundles on
an elliptic curve by Atiyah \cite{AtiyahElliptic},
there are two possibilities for the resulting rank two bundle. The first is
$\cO_E(p_1) \oplus \cO_E(p_1)$, and the second is the non-decomposable
rank two bundle $F_2$ defined by the
extension sequence
\be
0 \to \cO_E(p_1) \to F_2 \to \cO_E(p_1) \to 0
\ee
These two bundles define the same $S$-equivalence class,
but they are not isomorphic.
If we consider a generic family,
then we end up with $F_2$. It is called the regular
representative of the $S$-equivalence
class, meaning it is the representative with the
smallest automorphism group, or the
maximal Jordan block structure. Its Fourier-
Mukai transform, i.e. its $T$-dual, is $\cO_{2p_1}$,
the basic example of a degenerate brane with a gluing VEV.

Now we want to understand  the corresponding bundles on the
rational surface. This was studied by
\cite{FM} and \cite{ChenLeung}. Let us denote the line on $S_G$ corresponding
to $p_1$ by $l_1$ and the line corresponding
to $p_2$ by $l_2$. Then the corresponding rank two bundle
on $S_G$ is $\cO_S(l_1) \oplus \cO_S(l_2)$.
Now we let $p_1\to p_2$. As $l_1$ approaches $l_2$ and we get a
surface with an $A_1$ singularity.
We can blow up the surface to get an exceptional $-2$-curve $C$,
in the class $[l_2-l_1]$. Let us denote
the blow up by $\widetilde{S}_G$.

To extend the correspondence, we would like a tautological bundle
on $S_G$ (or $\widetilde{S}_G$) that naturally
arises as a limit of $\cO_S(l_1) \oplus \cO_S(l_2)$, has some kind
of interpretation as a regular representative,
and restricts to $F_2$ on $E$. One might think that one candidate
for such a rank two bundle
is $V = \cO_S(l_1) \oplus \cO_S(l_1 + C)$. However this clearly
restricts to $\cO_E(p_1)\oplus \cO_E(p_1)$ on
the boundary elliptic curve, not the regular representative. So
what is the tautological bundle that restricts to $F_2$?

Let us first consider $V$ restricted to the exceptional curve $C$.
Then we get that
\be
V|_C = \cO_C(-1) \oplus \cO_C(1).
\ee
This bundle is `unbalanced'. It admits a unique deformation to a
more generic rank two bundle
on $C$ with the same topological data, namely $\cO(0) \oplus \cO(0)$.
We can construct it explicitly
as an extension
\be
0 \to \ \cO_C(-1) \to \cO_C(0) \oplus \cO_C(0) \to \cO_C(1) \to 0
\ee
On a rational surface this deformation can be uniquely lifted to an
irreducible deformation
of the direct sum bundle $V$. To see this, we start with
the exact sequence
\be
0 \ \to \ \cO_{\widetilde S} \ \to \ \cO_{\widetilde S}(C) \ \to \ \cO_C(C) \ \to \ 0
\ee
The associated long exact sequence in cohomology gives us
\be
\ldots \to 0 \to \ H^1(\cO_{\widetilde{S}}(C)) \ \to \ H^1(\cO_C(C)) \to \ 0\ \to \ldots
\ee
where we used that $h^{1,0}(\widetilde S) = h^{2,0}(\widetilde S)= 0$.
But we have $\Ext^1(\cO(l_1),\cO(l_2)) \simeq H^1(\cO_{\widetilde S}(C))$ and
$\Ext_C^1(\cO(1),\cO(-1)) \simeq H^1(\cO_C(C))\simeq {\bf C}$, so our deformation on $C$ lifts
to $\widetilde S$, in fact we get the unique non-trivial extension
\be\label{Reg_extension_seq}
0 \to \cO(l_2) \to W \to \cO(l_1) \to 0
\ee
Note also that $\Ext^1(\cO(l_1),\cO(l_1)) = \Ext^1(\cO(l_2),\cO(l_2)) =0$ since
$h^{1,0}(\widetilde S)=0$,
so this is the only natural candidate for a bundle that restricts to $F_2$ on the boundary
and arises as a limit as $l_1 \to l_2$.
Therefore the deformed bundle $W$ is the tautological bundle
on $\widetilde{S}_G$ that corresponds
to the regular representative. Since it is trivial when restricted to $C$, it may also be thought
of as the pull-back of a rank two vector bundle on $S_G$ itself, as opposed to a more general sheaf.

For later purposes, we would like to note that
$\Ext^1_{\widetilde{S}}(\cO(l_1),\cO(l_2)) \cong H^1(\cO(l_2-l_1))$
vanishes when $l_2 - l_1$ is not effective. This is somewhat analogous to the
orthogonality of Fourier modes.

To see this, we can look at the index
\be
\sum_i (-1)^i \dim \Ext^i_{\widetilde{S}}(\cO(l_1),\cO(l_2))
\ee
which is invariant under deformations. Further,
$\Ext^2_{\widetilde{S}}(\cO(l_1),\cO(l_2))\cong H^0(\cO(l_1-l_2)\otimes K_{\widetilde{S}})$
vanishes on a rational surface, since $K_{\widetilde{S}}$
does not have sections and tensoring with $\cO(l_2-l_1)$ doesn't
make it any better. Therefore if we vary $l_1$ and $l_2$, $\Ext^1$ and $\Ext^0$ must appear or disappear in pairs.
But
$\Ext^0_{\widetilde{S}}(\cO(l_1),\cO(l_2))= H^0(\cO(C))$ is one-dimensional
when $l_2-l_1$ is effective, and zero dimensional when $l_2-l_1$ is not effective.
Therefore the unique $\Ext^1$ we found when $l_2-l_1$ is effective disappears when $l_2-l_1$ ceases to be
effective. An alternative proof which does not rely on deformations can be found in \cite{FM}.

\newpage

\newsection{Fibering the correspondence}

The correspondence discussed so far holds for individual rational surfaces.
Now we want to extend it to families, so that we can fiber this correspondence
over a base $B$.

\newsubsection{Constructing $Y_C$}

The first step is to construct a K\"ahler manifold
$Y \to B$ whose fibers are ADE rational surfaces.
We will assume that we are given
a smooth elliptically fibered K\"ahler manifold $Z$ with projection
$\pi_Z:Z \to B$ and section $\sigma$. We further assume
we are given a spectral cover
$C\subset Z$ which is a finite cover $C\to B$. We could
presumably generalize this set-up
by allowing the fibers of $Y \to B$ to be more general surfaces,
as long as their homology lattice contains an ADE root lattice,
and allowing for abstract unembedded covers $C \to B$.
In fact this would conceptually
be a cleaner way to describe the correspondence. But working with embedded covers
tends to make the set-up a little more concrete as we can write explicit equations.

To construct a bundle
we further need a spectral sheaf ${\cal L}$ on $C$, or a holomorphic bundle ${\cal V}$ on
$Z$ whose Fourier-Mukai yields such a spectral sheaf supported on $C$. However ${\cal L}$ plays no role
in constructing $Y$ itself, so it will be ignored in this subsection.

Let us first consider fibrations by surfaces of type $A_n$. Recall that we
constructed such surfaces by embedding an elliptic curve into ${\bf P}^2$
using the linear system $|3p_0|$ and then blowing up $n+1$ points $p_0, \ldots, p_n$.
The global analogue of this is as follows.
Consider the line bundle $\cO_Z(3\sigma)$ and its direct image $\pi_{Z*}\cO_Z(3\sigma)$.
This is a rank three bundle over $B$ with fibers $H^0(E_b, \cO_{E_b}(3p_0))$ over a
point $b \in B$. Given a point $p\in E_b$ and a section $s \in H^0(E_b, \cO_{E_b}(3p_0))$,
we have the natural evaluation map $(p,s)\to s(p)$, so dually we get
the embedding $E_b \hookrightarrow {\bf P}(H^0(E_b, \cO_{E_b}(3p_0))^\vee) \cong {\bf P}^2$.
Globally this means we get an embedding
\be
Z \ \hookrightarrow\ {\bf P}(\pi_{Z*}\cO_Z(3\sigma)^\vee)
\ee
where the right-hand-side is a ${\bf P}^2$-bundle over $B$.

Now $C$ sits in $Z$, and therefore by the above embedding
it also sits in ${\bf P}(\pi_{Z*}\cO_Z(3\sigma)^\vee)$.
Then $Y_C$ is obtained by blowing up
${\bf P}(\pi_{Z*}\cO_Z(3\sigma)^\vee)$ along $C\cup \sigma$:
\be\label{Y_C_def}
Y_C \ = \ {\rm Bl}_{C\cup \sigma}[{\bf P}(\pi_{Z*}\cO_Z(3\sigma)^\vee)]
\ee
This space is actually singular when
$C$ and $\sigma$ intersect. The reason is that along $C\cap \sigma$ we are effectively
blowing up along a length two scheme and such a blow-up yields a singular space.
To resolve this we can do a further blow-up of $Y_C$, or alternatively we can blow-up
along $C$ and $\sigma$ successively, which will do the extra blow-up for us (with exchanging
the order being related by a flop). The extra blow-up
is familiar from `matter curves' in $F$-theory compactifications.

We note that it is straightforward to do the blow-up explicitly. Our elliptic Calabi-Yau $Z$
is given by a cubic equation $c=0$ in the ${\bf P}^2$-bundle ${\bf P}(\cO\oplus L^2 \oplus L^3)$,
where $L = K_B^{-1}$ and
\be
c = -y^2 z + x^3 + f x z^2 + g z^3
\ee
The spectral cover $C$ in $Z$ is specified by an equation $s=0$ where
\be
s = p_0 + p_2 x + p_3 y + \ldots + p_{n+1} x^{(n+1)/2}
\ee
for $n$ odd, with a similar equation for $n$ even.
Here the $p_i$ are sections of $L^{n+1-i}$.
To perform the blow-up, we only need to introduce additional projective coordinates $(u,v)$ such
that $uc = vs$. The sections $p_i$ allow us to vary the rational surface fibration while keeping
the elliptic fibration fixed.

The Picard lattice of our surface fiber over a point $b\in B$ is of the form
\be
H^2(S_b, {\bf Z}) \ = \ \Lambda_{A_{n-1}} \oplus {\bf Z} E_b \oplus {\bf Z}l_0\oplus {\bf Z}f
\ee
and furthermore by construction these pieces do not mix as we vary the point $b$.

The construction for other Lie algebras is a little more involved. The cases $E_6$, $E_7$ and $E_8$ are
effectively discussed in \cite{Friedman:1997ih}. For $D_n$ we can give a fairly explicit description
using recent results on the Sen limit \cite{Clingher:2003ui}. Let $F \to B$ be ${\bf P}^1$-fibration
over $B$. Consider  ${\bf P}(\cO_F \oplus K_F \oplus L)$, which is a ${\bf P}^2$-fibration over $F$,
and let $(y,u,v)$ denote fiber coordinates of this ${\bf P}^2$-bundle. We take our rational surface fibration
$Y$ to be specified by an
equation
\be
y^2 = b_2 u^2 + 2 b_4 uv + b_6 v^2
\ee
where the $b_i$ are sections of $K_F^{-2}$,
$K_F^{-1}\otimes L^{-1}$ and $L^{-2}$ respectively. The line bundle $L = K_F^3$ is what appears in the usual Sen limit.
Geometrically the above equation cuts out a conic in each ${\bf P}^2$-fiber, and the fiber of
$Y \to B$ is a rational surface which is fibered by conics. For $v=0$ we recover the Calabi-Yau $Z$,
given by $y^2/u^2 = b_2$ (a double cover of $F$). The spectral cover $C$ is given by the intersection
of $Z$ with $R=\pi^{-1}(\Delta)$, where $\Delta = b_2 b_6 - b_4^2$ is a section of $K_F^{-2} \otimes L^{-2}$.
If $L$ has degree $k$ on the fiber of $F \to B$, then this cover has degree $4k+8$ over $B$.

In each case, the boundary of a rational surface fiber is the elliptic curve $E_b$,
so $Y_C$ contains an elliptic fibration over $B$ which is our original $Z$. The collection of lines on
each rational surface fit together in the cylinder $R\to B$, and the intersection of $R$ (after embedding
in $Y$) with $Z$ gives
the spectral cover $C \to B$.

Now given a spectral sheaf ${\cal L}$ supported on $C \subset Z$, associated to a bundle
${\cal V}$ on $Z$ by the usual Fourier-Mukai transform, we want to construct a natural
bundle ${\cal W}$ on $Y$. This bundle will have the property that restricting to the boundary gives us back ${\cal V}$,
i.e. ${\cal W}|_{Z} = {\cal V}$. The structure group of ${\cal W}$ will be a ${\bf C}^*$-extension
of the ADE structure group of ${\cal V}$.

In concrete applications we may encounter variants of the above situation. In one common situation,
some part of the gauge symmetry remains unbroken. Then the lattice $\Lambda_{rt}$ may further split up, eg.
$\Lambda_{rt} = \Lambda_{rt,1} \oplus \Lambda_{rt,2}$, and these pieces don't mix as we
vary over a base.  Then of course we could construct bundles associated
to the sublattices.

\newsubsection{General comments on the construction}

In order to understand the construction, it is useful to keep the conventional Fourier-Mukai
transform in mind, so let us briefly review this. Our discussion is somewhat heuristic,
in particular more attention should be paid to the branch locus of the spectral cover.
For a more detailed treatment see eg. \cite{Huyb_FM,BBR_FM}.

Let $Z\to B$
be a smooth elliptically fibered K\"ahler manifold with section $\sigma$.
Let us consider two copies of $Z$, denoted by $Z_1$ and $Z_2$. We want to think
of the elliptic fibers of $Z_1$ as being the dual of the elliptic fibers of $Z_2$.
More formally we can write
\be
Z_2\ =\ {\rm Pic}^0(Z_1/B)
\ee
where ${\rm Pic}^0(Z/B)$ is the relative Picard group of degree zero divisors. Since
${\rm Pic}_0(E) \cong E^\vee \cong E$ and $Z_1$ has a section, $Z_2$ is isomorphic to $Z_1$.
Now we construct the following double elliptic fibration
\be
Z_1 \times_{B} Z_2 \ = \ Z_1 \times_B {\rm Pic}_0(Z_1/B)
\ee
There are natural maps $\pi_1:Z_1\times_B Z_2 \to Z_1$ and $\pi_2:Z_1\times_B Z_2\to Z_2$ that
project on the two factors respectively.

On $Z_1\times_B Z_2$ we can define
a certain universal line bundle, the Poincar\'e line bundle ${\cal P}$, as follows:
we label a point in $Z_1\times_B Z_2$ as $(p,q)$, where $p\in Z_1, q\in Z_2$, and
$\pi_1(p) = \pi_2(q)$. We will also use $\bar p$ and $\bar q$ to denote
the corresponding points on the elliptic fibers $E = \pi_1^{-1}\pi_1(p)\cong
\pi_2^{-1}\pi_2(q)$. Since our elliptic fibrations have a section, there is a distinguished point
$\bar q_0$ on every elliptic fiber. Then the fiber of ${\cal P}$ at $(p,q)$ is given by
\be
{\cal P}_{(p,q)}\ =\ \cO_{E}(\bar q-\bar q_0)|_{\bar p}.
\ee
The Poincar\'e sheaf is the sheaf
of sections of ${\cal P}$, which we denote by the same name.

Now given a spectral sheaf ${\cal L}$ on $Z_2$, the dual bundle ${\cal V}$ on $Z_1$ is given
by
\be\label{FM_transform}
{\cal V} \ = \ \pi_{1*}(\pi_2^*{\cal L} \otimes {\cal P})
\ee
Note that if the spectral cover has degree $n$ over $B$, then $\pi_2^{-1}C= Z\times_B C$ has degree $n$ over $Z_1$, and
so $\pi_{1*}$ gives a rank $n$ bundle on $Z_1$.

We would like to make two comments. The first comment is that the Poincar\'e sheaf is defined on all of $Z_1\times_B Z_2$,
but in formula (\ref{FM_transform}) we actually only need the restriction of ${\cal P}$ to
$Z\times_B C \subset Z_1\times_B Z_2$, since we are tensoring with the sheaf $\pi_2^{*}{\cal L}$ which
is only supported on $Z\times_B C$. Therefore we may also formulate the Fourier-Mukai transform as follows.
Let us define
\be
\widehat{Z} \ = \ Z \times_B C
\ee
Due to the inclusion $Z\times_B C \to Z \times_B Z$ we may pull back the Poincar\'e sheaf to a sheaf
${\cal P}_{\widehat{Z}}$ on $\widehat{Z}$. Since the spectral sheaf is supported on $C$, we may then also
express ${\cal V}$ as
\be\label{FM_transform_II}
{\cal V} \ = \ {\bf FM}({\cal L}) \ = \ \pi_{1*}(\pi_2^*{\cal L} \otimes {\cal P}_{\widehat{Z}})
\ee
We will use this observation in the new construction.

The second comment we would like to make is that the Poincar\'e sheaf often fails to be a line bundle,
but this is not a cause for concern. Indeed the cross product $Z_1\times_B Z_2$ has singularities over points
in the base $B$
where we take the cross product of nodal elliptic curves, and the Poincar\'e sheaf fails to be a line
bundle at these singular points. We will see a very similar phenomenon in our new construction.

Now we would like to construct natural ADE bundles on a complex manifold $Y$
that admits a fibration $Y \to B$ by rational surfaces.
The idea of our construction will be to proceed along the same lines as above:
we define a suitable `tautological' sheaf
and then do an analogue of the Fourier-Mukai transform. Let us explain why such a tautological
sheaf should exist. We are grateful to T. Pantev for the argument below.

Consider the space $S\times {\rm Pic}(S)$. There is a universal sheaf ${\cal P}$
on this space:
given a point
$y\in {\rm Pic}(S)$, $y$ determines a sheaf ${\cal L}_y$ on $S$, and we take
\be
{\cal P}|_{S\times y}\ =\ {\cal L}_y.
\ee
We can easily extend this to families.
Given a fibration $Y \to B$, the relative Picard group ${\rm Pic}(Y/B)$ is defined as $R^1\pi_{*}\cO_Y^*$.
On a generic fiber $S_b$ this yields $H^1(S_b,\cO^*)$ which is the usual Picard group of $S_b$.
We define a new space by crossing with $Y$:
\be
{\cal Y}\ =\ Y_C \times_B {\rm Pic}(Y_C/B)
\ee
There is a universal sheaf ${\cal P}$ on ${\cal Y}$, as follows:
given a point $(x,y) \in {\cal Y}$ with $\pi(x)=\pi(y)=b$, $y$ determines a sheaf on
$S_b$, so we use this as the fibers of ${\cal P}$.

Now we define
\be
\widehat{Y}\ =\ Y_C \times_B C
\ee
It will be convenient to again use $\pi_1$ and $\pi_2$ for projection
on the first and second factor respectively.
There is a natural map from $C$ to ${\rm Pic}(Y_C/B)$, since by construction of $Y_C$
a point $p$ on
$C$  with $\pi_C(p)=b$ determines a line $l$ in $S_b$ and a generator $\cO_{S_b}(l)$
of ${\rm Pic}(S_b)$.
Then by pulling back the universal sheaf ${\cal P}$ we get a sheaf ${\cal P}_{\widehat Y}$ on
$\widehat{Y} = Y_C \times_B C$.

Once we have a Poicar\'e sheaf, we can proceed as in the usual Fourier-Mukai
transform, and define a bundle ${\cal W}$ on $Y_C$ as the transform
\be\label{ADE_transform}
{\cal W} \ = \ {\bf T}({\cal L}) \ =  \ \pi_{1*}(\pi_2^*{\cal L} \otimes {\cal P}_{\widehat{Y}})
\ee
Furthermore by construction, when we restrict this to $Z\times_B C \subset Y_C\times_B C$, we
reproduce the usual Fourier-Mukai transform for elliptically fibered manifolds reviewed above.

\newsubsection{Explicit description of the Poincar\'e sheaf}

\subseclabel{Poincare_def}

The construction thus boils down to defining a suitable generalized Poincar\'e sheaf ${\cal P}_{\widehat Y}$.
In this section we describe this Poincar\'e sheaf more explicitly.

Let us introduce the cylinder $R$, the collection of lines on each rational surface fibered over $B$.
$R$ has the structure
of a ${\bf P}^1$-fibration over $C$. After embedding in $Y$,
it can be identified with the exceptional divisor of the blow-up (\ref{Y_C_def}).
We also define
$\widehat{R} = R\times_B C$. Since the lines are naturally embedded in the surface fibers,
we have an embedding  $\widehat R \hookrightarrow \widehat{Y}$.
Furthermore since $R\to C$ is a ${\bf P}^1$-fibration,
we also get the ${\bf P}^1$-fibration
\be
\widehat{R} \ \to \ C\times_B C
\ee
Now we consider the diagonal $\Delta \in C\times_B C$. Pulling $\Delta$ back to $\widehat{R}$ and embedding
in $\widehat Y$ gives a natural divisor $D \subset \widehat Y$. The relations are summarized in
the following diagram:
\be
\begin{array}{cccccc}
   D&\hookrightarrow  & R \times_B C & \hookrightarrow & \widehat{Y} = Y \times_B C& \\
   \downarrow &  & \downarrow &  &  \downarrow & \!\!\!\!\!\searrow\\
   \Delta & \hookrightarrow & C \times_B C &  &  Y & C
\end{array}
\ee
Now our tentative definition
of the Poincar\'e sheaf is simply
\be
{\cal P}_{\widehat Y} \ = \ \cO_{\widehat Y}(D)
\ee

In practice one typically wants to consider a small modification of this, in order for the restriction of
the bundle ${\cal W}$ to $Z$ to have an ADE structure group instead of a ${\bf C}^*$-extension of this.
Let us do this explicitly for the the $A_{n-1}$ case. Fiber by fiber the correspondence
is essentially unchanged if we tensor $\cO_{\widehat Y}(D)$ by a line bundle pulled back from $C$ or from $Y$.
The first modification we want to make is for the Poincar\'e sheaf to be flat when restricted to the elliptic
fibers of $Z$.
We saw that this is achieved on the rational surface fibers
by tensoring with the line bundle $\cO(-l_0)$, which produces the collection
$\{\cO(l_1 - l_0), \ldots, \cO(l_n-l_0)\}$ instead
of $\{ \cO(l_1),\ldots, \cO(l_n)\}$. Fibering the line $l_0$ over $B$ gives the divisor $R_\sigma$ in $Y_C$.
Then we modify ${\cal P}$ as:
\be
{\cal P}_{\widehat Y} \ = \ \cO_{\widehat Y}(D - R_\sigma \times_B C)
\ee
Secondly, in order to produce the usual Poincar\'e sheaf on $Z\times_B C$,
we want to make this a little more symmetric.
We can do this by modifying ${\cal P}$ by the pull-back of
$\cO_C(\Sigma)$, where $\Sigma=\sigma\cap C$:
to
\be
{\cal P}_{\widehat Y} \ = \ \cO_{\widehat Y}(D - R_\sigma \times_B C - \pi_2^{-1}(\Sigma))
\ee
Note that we can think of $\pi_2^{-1}(\Sigma)$ as $(Y\times_B C) \cap (Y \times_B \sigma)$ in $Y\times_B Z$.
In this way we reproduce the more conventional expression $\cO(\Delta - \sigma\times_B C - \pi_2^{-1}\Sigma)$
when restricted to $Z\times_B C$. Finally, although our twist ensured that the restriction
${\cal W}|_Z$ has vanishing first Chern class along the elliptic fibers, it still has non-vanishing
first Chern class along the base $B$. To fix this we need to further twist by a factor
of $\pi_B^* K_B^{-1}$, where $\pi_B$ is the natural projection $Y \times_B C \to B$. Thus the final
expression that we settle on is given by
\be\label{P_sheaf_def}
{\cal P}_{\widehat Y} \ = \ \cO_{\widehat Y}(D - R_\sigma \times_B C - \pi_2^{-1}\Sigma)\otimes \pi_B^*K_B^{-1}
\ee
Most of the interesting behaviour however has to do with the factor $\cO_{\widehat Y}(D)$.

\newsubsection{Behaviour of the Poincar\'e sheaf at the branch locus}

There wasn't much choice in defining our Poincar\'e sheaf, the whole construction being essentially
tautological. It is not hard to see that our transform reproduces the correspondence
discussed in section \ref{Correspondence} for generic fibers. In that sense, the whole point now is to check
that we get the expected
behaviour also in more degenerate situations. In particular if we have a non-trivial fibration
$Y_C \to B$, then some of the fibers are going to be singular even in the generic case. This corresponds
to the branch locus of the spectral cover $C$, or equivalently the locus where two lines $l_1$ and $l_2$
on the surface fibers of $Y_C$ come together and get exchanged (which we may refer to as the branch locus of $Y_C$).

We want to take a closer look at the divisor $D$ in such a situation. It turns out that $\widehat Y$ is singular in codimension
three even for generic spectral covers, and $D$ is a Weil but not a Cartier divisor on $\widehat{Y}$. Thus
our ${\cal P}_{\widehat Y}$ is a rank one sheaf and not a line bundle, although it looks like a line
bundle away from the singular locus. As we observed previously, a similar phenomenon occurs for the conventional
Poincar\'e sheaf on $Z_1\times_B Z_2$.

To see this more explicitly, let's look at some local models. Let us denote by $t$ a local parameter on the base $B$.
Then we have the following local models for $C$ near the branch locus and for $Y$ when two lines come together:
\be
C: \ {\bf C}[t,s]/(t-s^2), \qquad Y: \  {\bf C}[t,x,y,z]/(t-x^2+y^2-z^2)
\ee
We then find the following local model for $\widehat Y = Y\times_B C$:
\be
\widehat{Y}: \ {\bf C}[t,s,x,y,z]/(t-s^2, t-x^2+y^2-z^2)\ \cong\ {\bf C}[s,x,y,z]/(s^2-x^2+y^2-z^2)
\ee
We see that $\widehat Y$ has conifold singularities in codimension three over the branch locus in $B$.

Now we can write $D$ explicitly in this local model for $\widehat Y$. We first write the corresponding
local models for $\widehat{R}$ and $C\times_B C$:
\be
\widehat{R}:\ {\bf C}[v,s_1,s_2]/(s_1^2 - s_2^2), \qquad C\times_B C: {\bf C}[s_1,s_2]/(s_1^2 - s_2^2)
\ee
Here we used the local model ${\bf C}[v,s_1,t]/(s_1^2 - t)$ for
$R$, with $v$ denoting a coordinate along the
${\bf P}^1$ fiber of $R \to C$. For a fixed $t$, this describes
the lines $l_1$ and $l_2$ that we met previously in section \ref{Correspondence}.

Now the diagonal $\Delta \subset C\times_B C$ is given by
\be
\Delta: \ {\bf C}[s_1,s_2]/(s_1-s_2)
\ee
This lifts to the subset
\be
D: \ {\bf C}[v,s_1,s_2]/(s_1-s_2) \ \subset \ \widehat{R}
\ee
Then finally we want to embed $D$ into $\widehat Y$. Equivalently
we find a map from the coordinate
ring of $\widehat Y$ to the coordinate ring of $D$, and we can take
this to be
\be
(x,y,z,s) \ \to (v,v,s_1,s_2)
\ee
Note that even locally we need not one but two equations to specify $D$ in
$\widehat{Y}$,
namely $x=y$ and $z=s$, so $D$ is not a Cartier divisor. Another
way to see this is that $D$ itself is clearly smooth, even when
it intersects the singular locus of $\widehat Y$. This means
that $D$ cannot be a Cartier divisor, since a Cartier divisor
would have to be singular at this locus.
Indeed the above is pretty much the standard example
of a Weil divisor that is not Cartier.

The conifold singularity of $\widehat{Y}$ can be removed by a small
resolution. It is perhaps interesting to note that the Poincar\'e
sheaf lifts to an actual line bundle if we pick one of the two small resolutions.

\newsubsection{Structure of ADE transform near the branch locus}

The next point we want to understand is the behaviour of the bundle we
constructed near the branch locus. In particular we would like to establish
two properties. We would like to show that over the branch locus,
${\bf T}({\cal L})$ reproduces
the tautological bundle corresponding
to the regular representative discussed
in section \ref{Reg_representative}. We would further like
to show that ${\bf T}({\cal L})$ is locally free,
so that it defines a smooth vector bundle on $Y$.
In this subsection we will establish both of these properties. This is
slightly tricky since we saw that $\widehat Y$ is actually singular
and ${\cal P}_{\widehat Y}$ fails to be a line bundle at the singularities.

The story is (mostly) local, and so we first want to compute
$(\pi_{\widehat{Y}/Y})_{*}\cO_{\widehat Y}(D)$ in the local models given above.

To simplify the discussion, let us also consider the divisor $D'$ defined as
follows:
\be
D': \ {\bf C}[x,y,z,s]/(x-y,z+s)\quad \subset\quad \widehat{Y}
\ee
It is also a Weil divisor, but the sum $D+D'$ is defined by a single equation
(namely $x=y$) and so corresponds to a Cartier divisor. Thus although
neither $\cO(D) $ or $\cO(D')$ is a line bundle, $\cO(D+D')$ actually is a
line bundle on $\widehat{Y}$. Its sections are local meromorphic functions
that are allowed to have a pole at $x=y$, in other words expressions of type
\be
f(x,y,z,s)/(x-y)
\ee
where $f$ is holomorphic. Furthermore, $\cO(D+D')$ is actually the
pull-back of the line bundle $\cO_Y(R)$ downstairs, since
$\pi_{\widehat{Y}/Y}^{-1}(R) = D + D'$. Indeed the equation of $R$ in $Y$ is given simply by
$x=y$, describing a pair of lines $l_1$ and $l_2$ for each $t$.

Now we can start computing the push-forwards. We write
\be
\cO_{\widehat{Y}}(D) \ = \ \cO_{\widehat{Y}}(D+D') \otimes \cO_{\widehat{Y}}(-D')
\ee
so that
\be
\pi_{\widehat{Y}/Y*}\cO_{\widehat{Y}}(D)\ =\ \cO_Y(R) \otimes \pi_{\widehat{Y}/Y*}\cO(-D')
\ee
The advantage of writing it this way is that $\cO_{\widehat{Y}}(-D')$ is an ideal sheaf,
and the push-forward of an ideal sheaf is a bit easier to understand.

We will use the standard short exact sequence associated to a divisor:
\be
0 \ \to \ \cO_{\widehat{Y}}(-D') \ \to \cO_{\widehat{Y}} \ \to \ \cO_{D'} \ \to \ 0
\ee
In order to save some notation, we are going to denote the map
$ \pi_{\widehat{Y}/Y}$ simply by $\pi$ in the remainder of this subsection.
Applying the push-forward we get
\be
0 \ \to \ \pi_*\cO_{\widehat{Y}}(-D') \ \to \pi_*\cO_{\widehat{Y}} \ \to \ \pi_*\cO_{D'} \ \to \ 0
\ee
The sequence terminates because $\pi$ is a finite map, hence $R^1\pi_* = 0$.
Let us first understand $\pi_*\cO_{\widehat{Y}}$. Sections of $\cO_{\widehat{Y}}$ are of the form
\be
1\cdot f(x,y,z,s)
\ee
where $f$ is any element of the ring ${\bf C}[x,y,z,s]/(x^2-y^2+z^2-s^2)$. Applying
$\pi_*$ means that we should regard any such section as a generator over the
ring ${\bf C}[x,y,z,t]/(x^2-y^2+z^2-t)$ where $t=s^2$. Clearly we can write any such section
as a combination of two generators, $1$ and $s$. Furthermore any appearance of $t$ can be eliminated
using the ring relation (equivalently any occurrence of $s^2$ can be eliminated using the ring relation
even before pushing forward). Thus we can write sections of $\pi_*\cO_{\widehat{Y}}$ as
\be
f_1(x,y,z) + s f_2(x,y,z)
\ee
The map $\cO_{\widehat{Y}}  \to  \cO_{D'}$ simply imposes the relations $x-y=z+s=0$, so we may
use this to eliminate $y$ and $s$ altogether, i.e our section maps as
\be
f_1(x,y,z) + s f_2(x,y,z) \ \to \ (f_1(x,y,z) - zf_2(x,y,z))|_{y=x}
\ee
The kernel of this map consists of sections $f_1(x,y,z) + s f_2(x,y,z)$ where
\be
f_1(x,y,z)\ =\ z f_2(x,y,z) + (x-y) f_3(x,y,z)
\ee
Rearranging as
\be\label{Reg_sections}
f_1 + s f_2\ =\ (z+s)f_2 + (x-y) f_3
\ee
we see that $\pi_*\cO_{\widehat{Y}}(-D')$ is freely generated over ${\bf C}[x,y,z]$
or ${\bf C}[x,y,z,t]$ by $x-y$ and $z+s$
(recall $t$ can always be eliminated using the ring relation). This is probably
not so surprising
given that $D'$ is defined on $\widehat{Y}$ as $x-y=z+s=0$,
although we have to be careful to say
which ring we are working over. At any rate, we see that
in our local model $\pi_*\cO_{\widehat{Y}}(-D')$ is a
locally free sheaf (a vector bundle) of rank two. In fact from (\ref{Reg_sections})
and the fact that the push-forward of the
Poincar\'e sheaf is given by tensoring
$\pi_*\cO_{\widehat{Y}}(-D')$ with an honest line bundle,
we see that locally we have the isomorphism
\be\label{P_loc}
(\pi_{\widehat{Y}/Y})_{*}{\cal P}_{\widehat{Y}}\ \cong\ \cO_Y \oplus \cO_Y
\ee
This is the first property we set out
to prove in this subsection.

In order to establish the second property, let us restrict our bundle to the surface
fibers of the fibration $Y_C \to B$. It is straightforward to see that
away from the branch locus, the restriction yields a direct sum of line bundles.
We need to show that when the bundle
we constructed is restricted to a surface fiber over the branch locus (here $t=0$ in the local model),
it can be expressed as the non-trivial extension in equation (\ref{Reg_extension_seq}).
The argument has a global and a local component. We start the local one.

Let us take a closer look at the bundle $\pi_{\widehat{Y}/Y*}\cO(-D')$
restricted to the surface fiber over $t=0$.
In order to save some notation, let us simply
denote $ \pi_{\widehat{Y}/Y*}\left.\cO(-D')\right|_{t=0}$ by $F_2$. We also denote the
fiber over $t=0$ (with equation $x^2-y^2+z^2=0$) by $S_0$, and its resolution by
$\nu:\widetilde{S}_0\to S_0$.
We now
claim that when we lift $F_2$ to $\widetilde{S}_0$,
it can be expressed as the extension of two specific line bundles $\cO_{\widetilde{S}_0}(-l_1)$
and $\cO_{\widetilde{S}_0}(-l_2)$,
thereby recovering the `regular representative' discussed in section \ref{Reg_representative}.

To see this, let us take general sections of $F_2$ as expressed in (\ref{Reg_sections}),
and consider the restriction mapping $s \to 0$.
Note this is consistent with $s^2 = 0$ at $t=0$. We are left with
\be
 z\, f_2 + (x-y) f_3
\ee
which are precisely the sections of the ideal sheaf ${\cal I}_{(x-y,z)}$ on $x^2 - y^2 + z^2 = 0$.
Thus the map $s\to 0$ describes an onto map from $F_2$ to ${\cal I}_{(x-y,z)}$. We can express
this as a sequence:
\be
F_2\ \mathop{\to}^r \ {\cal I}_{(x-y,z)}\ \to\ 0.
\ee
We would like to complete this to a short exact sequence:
\be
0 \ \to \ {\cal K} \ \to \ F_2\ \mathop{\to}^r \ {\cal I}_{(x-y,z)}\ \to\ 0.
\ee
which we will eventually claim is essentially the extension sequence
(\ref{Reg_extension_seq}) we are trying to recover.
To do this we need to look
at the kernel of the map $r$. Again we identify this sheaf by its sections. The kernel is a sheaf with sections
of the form (\ref{Reg_sections}) which are in the kernel of the map $s\to 0$. Such sections
are of the form
\be
s\,f_2(x,y,z)
\ee
where $z\, f_2 + (x-y) f_3=0$. This defines ${\cal K}$, albeit somewhat implicitly.

In order to get a more recognizable description of the kernel sheaf, we will lift to the blow-up
$\widetilde{S}_0$,
where we can solve the equation $z\, f_2 + (x-y) f_3=0$ more easily.
To do this we introduce projective coordinates $(\lambda_1,\lambda_2)$, and instead of $x^2 - y^2 + z^2 =0$
we impose the pair of equations
\be\label{Szero_Resolution}
\left(
  \begin{array}{cc}
    x-y & z \\
    -z & x+y \\
  \end{array}
\right)
\left(
  \begin{array}{c}
    \lambda_1 \\
    \lambda_2 \\
  \end{array}
\right) = 0
\ee
on ${\bf C}^3 \times {\bf P}^1$. On the blow-up we can solve $z\, f_2 + (x-y) f_3=0$ by writing
\be
f_3\ =\ \lambda_1\, q(x,y,z,\lambda), \qquad f_2\ =\ \lambda_2\, q(x,y,z,\lambda),
\ee
where $q$ is an arbitrary polynomial. Then the kernel
is isomorphic to a sheaf with sections
\be
\lambda_2\, q(x,y,z,\lambda)
\ee
which defines the ideal sheaf
${\cal I}_{(\lambda_2)}$, and the map
from ${\cal I}_{(\lambda_2)}$ to $F_2$ is given by multiplying such
an expression by $s$, i.e. it is the inclusion
\be
i: \lambda_2 q \ \to \ s(\lambda_2 q )
\ee
Thus we have derived the short exact sequence
\be
0 \ \to \ {\cal I}_{(\lambda_2)} \ \mathop{\to}^i \ \nu^*F_2 \ \mathop{\to}^r \ \nu^*{\cal I}_{(x-y,z)} \ \to \ 0
\ee
But as one may check from equation (\ref{Szero_Resolution}), $\lambda_2 = 0$ also implies that
$x-y = z = 0$, so $\lambda_2=0$ is the equation for the line
for the line $l_1$
on $\widetilde{S}_0$ (the proper transform of the line
$x-y=z=0$ downstairs on ${S}_0$).
Thus the ideal sheaf
${\cal I}_{(\lambda_2)}$ is in fact $\cO(-l_1)$. Therefore when pulled back to the resolution of
$x^2-y^2+z^2=0$, we found the following short exact sequence:
\be
0 \ \to \ \cO(-l_1) \ \mathop{\to}^i \ \nu^*F_2 \ \mathop{\to}^r \ \nu^*{\cal I}_{(x-y,z)} \ \to \ 0
\ee

Furthermore, when we pull back
${\cal I}_{(x-y,z)}$ to the resolution, we are looking at $\cO(-l_2)$,
where $l_2$ is the locus $x-y=z=0$ on $\widetilde{S}_0$, equivalently
the pull-back of the line $x-y=z=0$ downstairs on ${S}_0$. Thus $l_2$ corresponds
to $l_1 \cup C$, where $C$ is the exceptional ${\bf P}^1$ of the blow-up.
Putting this together, we see that after lifting to the blow-up, we have found the
short exact sequence
\be
0 \ \to \ \cO_{\widetilde{S}_0}(-l_1) \ \to \ \nu^*F_2 \ \to \ \cO_{\widetilde{S}_0}(-l_2) \ \to \ 0
\ee
with $l_2 = l_1 + C$.

We now want to argue that the extension is non-trivial, i.e. $\nu^*F_2$ is not
isomorphic to the direct sum
$\cO(-l_1) \oplus \cO(-l_2)$. For this we really need to compactify to a rational surface,
since higher cohomology vanishes on affine varieties and the sheaf would just be isomorphic to $\cO\oplus \cO$
in the local model.
The precise compactification doesn't
matter much, as long as $l_1$ and $l_2$ are promoted to $-1$-curves on a rational surface,
because then we can use the argument in
section \ref{Reg_representative}, which says that it is enough to check the extension over the
exceptional curve $C$. For example a simple way to compactify the surface fibers of our local model is as follows:
\be\label{C_Bundle}
x^2 - y^2 + w^2 (z_1^2 - t z_0^2) = 0
\ee
Here $(z_0,z_1)$ are projective coordinates on ${\bf P}^1$ and $(x,y,w)$ are projective coordinates on
${\bf P}^2$. For each $t$ the equation (\ref{C_Bundle}) defines a conic bundle over ${\bf P}^1$, and in the patch
$z_0 = w = 1$ we recover our previous local model. This can be viewed as a Hirzebruch surface blown up at two points,
so it fits with the $A_n$ and $D_n$ rational surfaces discussed previously.
The fiber over a generic point $(z_0,z_1)$ is a smooth conic,
which degenerates to a pair of lines over $z = \pm \sqrt{t}$. Denoting the generic fiber by $f$, we may denote the
two lines over $z = \sqrt{t}$ by $l_1$ and $f-l_1$ respectively, since their sum deforms to a generic fiber. It follows
that $0 = f^2 = ((f-l_1)+l_1)^2$, and using the symmetry between $l_1$ and $f-l_1$ and the fact
that they intersect at a single point, we deduce that $l_1^2 = -1$. The same argument applies at $z=-\sqrt{t}$,
so the lines $x=y$, $z = \pm \sqrt{t}$ are promoted to $-1$-curves. They
collide on the fiber over $t=0$, creating the $A_1$ singularity.

To argue that the extension is non-trivial over $C$ is straightforward.
Using the fact that $l_1$ and $l_2$ have self-intersection $-1$, we have
$\cO(-l_1)|_C = \cO_C(-1)$ and $\cO(-l_2)|_C = \cO_C(1)$. Furthermore we
have already shown, $\nu^*F_2$ is the pull-back
of a vector bundle downstairs (a locally free sheaf) and so it restrict to the trivial rank two bundle
on $C$. Therefore on $C$ we get the short exact sequence
\be
0 \ \to \ \cO_C(-1) \ \to \ \cO_C(0)\oplus \cO_C(0) \ \to \ \cO_C(1) \ \to \ 0
\ee
which indeed describes the non-trivial extension of $\cO_C(1)$ by $\cO_C(-1)$, as required.

We are now essentially done. We only need to remember that $\nu^*F_2$ corresponds to
$\pi_{\widehat{Y}/Y*}\cO(-D')$ at $t=0$, but the push-forward
of our Poincar\'e sheaf is $\pi_{\widehat{Y}/Y*}\cO(D)$.
To save on notation, we denote the latter by $\widetilde{F}_2$.
The
difference is given by tensoring with the line bundle $\cO(R)$,
i.e. we have
\be
(\pi_{\widehat{Y}/Y})_{*}{\cal P}_{\widehat{Y}}\ =\ \widetilde{F}_2\ =\ \cO_{S_0}(R) \otimes F_2.
\ee
Here $R$ is the divisor defined by the equation $x=y$, which describes the locus
$x-y=z^2=0$ in the fiber over $t=0$, which in turn describes
the divisor $2l_1 + C = l_1 + l_2$ on the resolution $\widetilde{S}_0$. Therefore
after tensoring with $\cO(l_1 + l_2)$ we have
\be\label{P_ext}
0 \ \to \ \cO_{\widetilde{S}_0}(l_2) \ \to \ \nu^*\widetilde{F}_2 \ \to \ \cO_{\widetilde{S}_0}(l_1) \ \to \ 0
\ee
This is precisely the non-trivial extension sequence (\ref{Reg_extension_seq})
that we wished to recover,
so our transform has the expected behaviour at the branch locus. This concludes the second
property that we wanted to show.

\medskip

\noindent{\it Acknowledgements:} We are grateful to T.~ Pantev for discussion.
RD acknowledges partial support by NSF grant DMS 1304962. The research of MW
was supported in part by a visiting faculty position at the mathematics department of the University
of Pennsylvania.

\end{document}